\documentclass{amsart}
\usepackage{amsmath,amssymb,amsthm,latexsym}
\title{A hyperbolic surface with a square grid net}%

\usepackage{graphicx}
  \ifx\LabelFigloaded\MYundefined\relax
  \else
   \fi

  \def\LabelFigloaded{\relax}% now loaded

  \chardef\LabelFigCatAt\the\catcode`\@
  \catcode`\@=11

 \let\LabelFigwlog@ld\wlog
 \def\wlog#1{\relax}

 \ifx\\\MYundefined@
    \let\\\relax
 \fi

  \def\ms@g{\immediate\write16}

 \def\N@wif{\csname newif\endcsname }
 \def\Temp@ {\N@wif\ifIN@}
 \ifx\INN@\MYundefined@
    \else \let\Temp@\relax
 \fi
 \Temp@

  \def\IN@{\expandafter\INN@\expandafter}
  \long\def\INN@0#1@#2@{\long\def\NI@##1#1##2##3\ENDNI@
    {\ifx\m@rker##2\IN@false\else\IN@true\fi}%
     \expandafter\NI@#2@@#1\m@rker\ENDNI@}
  \def\m@rker{\m@@rker}
 
  \newtoks\Initialtoks@  \newtoks\Terminaltoks@
  \def\SPLIT@{\expandafter\SPLITT@\expandafter}
  \def\SPLITT@0#1@#2@{\def\TTILPS@##1#1##2@{%
     \Initialtoks@{##1}\Terminaltoks@{##2}}\expandafter\TTILPS@#2@}

 \def\Shifted@@#1#2#3{\setbox0=\hbox{#3}%
   \raise -\dp0\vbox {\kern-#2%
       \hbox {\kern#1\unhbox0\kern-#1}%
           \kern#2}}

 \newcount\gridcount
 \newbox\auxGridbox@ \newbox\hGridbox@ \newbox\vGridbox@
 \newbox\Labelbox@ \newbox\auxLabelbox@
 \newbox\Coordinatebox@
 \newtoks\Labeltoks@
 \newdimen\Wdd@ \newdimen\Htt@
 \newdimen\Wddd@ \newdimen\Httt@
 
 \def\Wr@{\immediate\write16}

 \newdimen\GL@wd%% grid-line width
 \def\GridLineWidth#1{\GL@wd=#1}

 \def\gobble#1{}
 \def\EdgeErr@{\Wr@{}%
      \Wr@{\string\Edges\space argument
      1, 10, 100 or 1000 please\string!}%
      }

 \newcount\Edgect@

 \def\Sweepup#1\endSweepup{}

 \def\SetEdges@{%
    \edef\Zr@@s{\expandafter\gobble\number\Edgect@\empty}%
        \count255=0\Zr@@s\relax
        \ifnum\count255=\z@\else\EdgeErr@\show\tailtest\fi
        \count255=1\Zr@@s\relax%\showthe\count255
        \ifnum\count255=\Edgect@\relax\else\EdgeErr@\show\leadtest\fi
    \EdgGl@b\edef\Zr@s{\expandafter\gobble\Zr@@s\empty}%\show\Zr@s
    \ifnum\Edgect@>\@ne\relax\EdgGl@b\let\L@Dc\empty
        \else\EdgGl@b\edef\L@Dc{\string.}\fi
    \ifnum\Edgect@>\@ne\relax
        \EdgGl@b\edef\Edgescale@##1{\divide##1 by \Edgect@}%
        \else\EdgGl@b\edef\Edgescale@##1{}\fi
    }

 \def\Edges#1{\Edgect@=#1\relax
     \let\EdgGl@b\global \SetEdges@}

 \Edges{1}%% default

 \def\hhrule{\hrule height \GL@wd\vskip-.\GL@wd}

 \def\hRule@{%
   \advance\gridcount -2%
   \vfil\hhrule\vfil
   \llap{\smash{\raise -2.5pt
     \hbox{\L@Dc\number\gridcount\Zr@s\kern2pt}}}%
   \hhrule
   }

\def\vvrule{\vrule width \GL@wd \kern-\GL@wd}

 \def\vRule@{\advance\gridcount 2%
   \hfil\vvrule\hfil
   \setbox\auxGridbox@=\vbox to 0pt
      {\vskip \Htt@\vskip 2pt
        \hbox to 0pt{\hss\L@Dc\number\gridcount\Zr@s\hss}\vss}%
      \wd\auxGridbox@=0pt \box\auxGridbox@
   \vvrule
   }

 \def\PlaceGrid@@{\gridcount=10 
  \setbox\hGridbox@=\hbox{%
        \hbox{%
             \hskip-.4pt\vrule
             \vbox to \Htt@{%
               \offinterlineskip\parindent=\z@\relax
               \hbox to \Wdd@{\hfil}
               \hRule@\hRule@\hRule@\hRule@
               \vfil\hhrule\vfil}%
             \vrule\hskip-.4pt}
    }%
  \gridcount=0%
  \setbox\vGridbox@=\hbox{%
      \vbox{\offinterlineskip\parindent=0pt\hsize=0pt
         \vskip-.4pt\hrule%
         \hbox to \Wdd@{%
                 \vtop to \Htt@{\vfil}%
                 \vRule@\vRule@\vRule@\vRule@
                 \hfil\vvrule\hfil}%
         \hrule\vskip-.4pt}}%
  \wd\hGridbox@=0pt\ht\hGridbox@=0pt
  \wd\vGridbox@=0pt\ht\vGridbox@=0pt
  \hbox{\box\hGridbox@\box\vGridbox@}%
  }

 \def\LabelsGlobal{\def\LabGl@b{\global}}
 \def\LabelsLocal{\def\LabGl@b{}}
 \LabelsGlobal %% default

 \def\SetLabels#1\endSetLabels{%
   \LabGl@b\Labeltoks@={#1()\\}%
   }

 \LabGl@b\Labeltoks@={()\\}

 \def\ShowGrid{\LabGl@b\let\PlaceGrid@\PlaceGrid@@}
 \def\HideGrid{\LabGl@b\let\PlaceGrid@\relax}
 \def\Grids{\ShowGrid\LabGl@b\let\GridSwitch@\ShowGrid}
 \def\noGrids{\HideGrid\LabGl@b\let\GridSwitch@\HideGrid}

 \noGrids

 \def\bAdjust@@{%
     \setbox\auxLabelbox@=\hbox{\raise \dp\auxLabelbox@
            \box\auxLabelbox@}}
 \def\bAdjust@{\let\vAdjust@\bAdjust@@}

 \def\eAdjust@@{\dimen0=-.5\ht\auxLabelbox@
     \advance\dimen0 by .5\dp\auxLabelbox@
     \setbox\auxLabelbox@=
            \hbox{\raise\dimen0\box\auxLabelbox@}}
 \def\eAdjust@{\let\vAdjust@\eAdjust@@}

 \def\tAdjust@@{%
     \setbox\auxLabelbox@=\hbox{\raise-\ht\auxLabelbox@
            \box\auxLabelbox@}}
 \def\tAdjust@{\let\vAdjust@\tAdjust@@}

 \let\vAdjust@\relax

 \def\lAdjust@{\let\hAdjust@\rlap}
 \def\rAdjust@{\let\hAdjust@\llap}

 \let\hAdjust@\relax\let\vAdjust@\relax

 \def\FetchLabel@#1(#2)#3\\{%
     \IN@0#2@@\ifIN@
        \setbox0=\hbox{\ignorespaces#1#3\unskip}%
        \ifdim\wd0>0pt
           \ms@g{}%
           \ms@g{ !!! Bad label(s)? !!!}%
           \message{ #1(#2)#3}%
        \fi
        \def\LabelMole@##1\endFetchLabel@{%
            \IN@0()\\@##1@%
            \ifIN@\def\Temp@{\FetchLabel@##1\endFetchLabel@}%
            \else\def\Temp@{}%
            \fi
            \Temp@
           }%
     \else
       \ignorespaces#1\unskip
       \setbox\auxLabelbox@=%
         \hbox to 0pt{\hss\ignorespaces\hAdjust@
          {\ignorespaces#3\unskip}\hss}%
       \vAdjust@
       \let\hAdjust@\relax\let\vAdjust@\relax
       \AugmentLabelBox@@{#2}%
       \ht\Labelbox@=0pt\dp\Labelbox@=0pt
       \let\LabelMole@\FetchLabel@%
     \fi\LabelMole@}

 \newtoks\XYSep@ %\XYSep@{*}
 \def\SetXYSeparator#1{%
     \IN@0#1@@\ifIN@\XYSep@{*}%
     \else
     \XYSep@{#1}%
     \fi
     }

 \SetXYSeparator*

 \def\AugmentLabelBox@@#1{%
     \IN@0\the\XYSep@ @#1@\ifIN@
       \SPLIT@0\the\XYSep@ @#1@%
       \setbox\Labelbox@=\hbox to 0pt{%
         \unhbox\Labelbox@
         \Shifted@@{\the\Initialtoks@\Wddd@}%
         {\the\Terminaltoks@\Httt@}%
         {\box\auxLabelbox@}}%
     \else
         \ms@g{}%
         \ms@g{ !!! Bad insertion point. !!!}%
         \message{ (#1\ this point was rejected.)}%
     \fi
    }

 \def\FetchOption@#1[#2]#3\endFetchOption@{%
    \def\temp{#1}%\show\temp
    \ifx\temp\empty
       \Edgect@=#2\relax%\showthe\Edgect@
       \let\EdgGl@b\relax
       \SetEdges@%\def\Edgescale@##1{\divide##1 by \Edgect@\relax}%
       \Cleaner@#3%
    \fi}

 \def\Cleaner@#1[@]{\Labeltoks@{#1}}
     
 \def\PlaceLabels@@{\mathsurround=0pt%\bgroup
     \def\Cr@{\\}%
     \let\L\lAdjust@\let\R\rAdjust@
     \let\B\bAdjust@\let\E\eAdjust@\let\T\tAdjust@
     \expandafter\FetchOption@\the\Labeltoks@[@]\endFetchOption@
     \Wddd@=\Wdd@ \Edgescale@\Wddd@ %\showthe\Edgect@
     \Httt@=\Htt@ \Edgescale@\Httt@
     \expandafter\FetchLabel@\the\Labeltoks@\endFetchLabel@
     \box\Labelbox@%\egroup
     }%

 \let \PlaceLabels@\PlaceLabels@@

 \def\AffixLabels#1{\setbox\Coordinatebox@=\hbox{#1}%
      \Wdd@=\wd\Coordinatebox@ \Htt@=\ht\Coordinatebox@
      \advance\Htt@ \dp\Coordinatebox@
      \hbox{\copy\Coordinatebox@\kern-\Wdd@ 
           \Shifted@@{0pt}{-\dp\Coordinatebox@}%
           {\PlaceLabels@\PlaceGrid@}%
           \kern\Wdd@}%
      \GridSwitch@ %% next grid hidden
      \LabGl@b\Labeltoks@{()\\}%
      }
 
   \let\wlog\LabelFigwlog@ld   %%restore logging
   \catcode`\@=\LabelFigCatAt  %%12 or 13

                                By

              Raymond S\'eroul <A18645@FRCCSC21.BITNET>
                                and 
              Laurent Siebenmann <lcs@topo.math.u-psud.fr>
    
              VERSIONS: July 1991, Oct 1991, Jan 1992, July 1992

INTRODUCTION

      This labelling package is intended for TeX users who
rely on non-TeX sources for for their graphics inserts.  It
provides means for adding TeX labels to such inserts with a
minimum of fuss. 

       For most labels, TeX users have in the past found it
reasonably convenient to rely on non-TeX sources. Typical
occasions when an inescapable need for TeX labels seemed to
arise are

 (a) when the graphics program lacks certain exotic or complex
mathematical symbols

 (b) when the very highest typographical quality is wanted for the
labels

 (c) when labels included with the graphics fail to print, 
 and you cannot figure out why (cf. boxedeps.doc).  The labels
 provided by labelfig.tex are 100% portable.

       Since this package first appeared, many users, who in the
past scarcely dreamed of using TeX labels, have come to use
nothing but.  So it is now appropriate to add

Intoxication Warning:  TeX labels may be addictive and expensive. 

     If you have a fast preview you may disagree, and even find
that this package provides an agreeable paste-up environment; see
extra applications at end.

     Note to publishers: It is possible and convenient to ultimately
export the TeX labels produced by labelfig.tex to become an integral
part of the EPS file. This is often desired by a publisher who typically
uses an "upmarket" graphics or page layout program, with which the
staff is skilled in perfecting figures.  See Appendix I for
a recipe.

     The authors are grateful to Patrick Ion of Math Reviews for
helpful comments and encouragement.

BASIC INSTRUCTIONS

    After reading in the macro file using

preview or proof your figure with a coordinate grid printed on
top, by typing the following:

    \ShowGrid  % shows grid  for next figure only
    \AffixLabels{<the graphics insertion>}

Here <the graphics insertion> is what you would type to insert
the graphics object alone without the grid.  This must provide
for the space around it. For example <the graphics insertion>
might well be \BoxedEPSF{MyFigure scaled 700} using the
boxedeps.tex macro package (from same source); this provides a
TeX box containing the encapsulated PostScript insert specified by
the file MyFigure. \AffixLabels{...} provides the grid (supposing
\ShowGrid is present) and later, once you have specified labels
using the grid, it will "tack on" the labels.

     The grid is a sort of (usually elongated) checkerboard of
ten rows and ten columns and its (internal) partitions are by
default numbered  .1, ... ,.9  both horizontally (X-coordinate
running left to right) and vertically (Y-coordinate running bottom
to top).  Thus the points enclosed by the grid correspond to the
points of the unit square in the cartesian "X-Y" plane, the lower
left corner corresponding to the origin (0,0).  By extrapolation,
the full page corresponds to a larger rectangle in the plane.

     These coordinates serve to position labels as follows.
Before the \AffixLabels{...} command type label specifications:

  \SetLabels
   (<X-coordinate>*<Y-coordinate>) <first label> \\
   .
   .
   .
   (<X-coordinate>*<Y-coordinate>)  <last label> \\
  \endSetLabels

Each row specifies one label and is terminated by \\.  In each
row, the position indicator comes first; it is written as a
standard cartesian point except that the X- and Y- coordinates
are separated by * rather than a comma because TeX allows a
comma as decimal point. There are no dimension units to specify
as the unit is the grid itself.

     By default, this cartesian point specifies where the middle
of the baseline of the label will be located.  However if you precede
the point by \L [or \R] the left [or right] edge of the baseline will
be located there. Similarly you may also precede the point by \T, \E,
or \B to vertically align the top equator or bottom of the label box
at the specified point.  This gives nine standard positions of
the label with respect to the insertion point --- corresponding to
the eight principle points of the compas and the center

                     \L\T     \T      \R\T

                     \L\E     \E      \R\E

                     \L\B     \B      \R\B

But this neglects the default "baseline" level of TeX,
giving potentially three more positions

                     \L    <no tag>   \R

For text, the baseline level is often the preferred. Its relation to
the others is variable. It will often coincide with the bottom level,
as happens for "X".  But it is often distinct, as for "g", in which
case you have in all 12 distinct positions rather than 9.

     It is convenient to think of this specification of label
position as attaching the label by a thumb-tack to the coordinate
grid. There are up to twelve positions of the thumb-tack on the
label, while the position of the thumb-tack on the coordinate grid is
arbitrary.  Normally, one choses the position of the thumb-tack on
the label to be the one that is the closest to the item being
labeled.  There are good reasons for this "rule of thumb":

   (a)  It facilitates correct positioning at first try.

   (b)  If the scale of the figure must be altered after labels
have been affixed, the labels have a good chance of remaining well
positioned.

   (c)  The visible grid need not extend beyond the "bounding box"
for the figure, because the best preferred position is always
(at least almost) within the bounding box .

The second reason is particularly important. Indeed it often
happens that scale has to be altered after labelling begins, in
order to either provide space for the labels, or to adjust
proportions between the labels and the figure.  (The size of labels
is unaffected by scaling.)

     Here is an artificial but self-contained test which uses
TeX rules to make a graphics object.

TEST

    Do not skip this!

 \def\FrameIt#1{\hbox{\vrule$\vcenter {\hrule\kern3pt%
             \hbox {\kern3pt #1\kern3pt}%
               \kern3pt\hrule}$\relax\vrule}}

 \def\Caption#1#2{\FrameIt{%
       \vtop {\hsize=#1\relax \parindent=0pt
         \leftskip=0pt \rightskip=0pt plus15pt
         \parfillskip=0pt
         \lineskip=1pt\baselineskip=0pt
         #2}}}

 \def\FirstQuadrant{\hbox to 100pt{\vrule\vbox to 100pt{%
        \hbox to 100pt{\hfil}\vfil\hrule}\hss}}

  \SetLabels
    \R(.5*.2) $\zeta\,\cdot$\\
    (.9*-.10) $\xi$\\
    \R(-.03*.9) $\eta$\\
    \T(.5*.9) \Caption{70pt}{%
          \it The norm of
          $g(\xi+i\eta)$ is indicated on
          contours of this invisible surface.}\\
  \endSetLabels

  \AffixLabels{\FirstQuadrant} \end

  Note that the coordinates to use for labels are indicated on the
edges of the grid (when visible) corresponding to the conventional
x- and y- axes of the Cartesian plane. By default the grid is
1-by-1. However, by the command \Edges{100}, you can change this
to 100-by-100 and many users find this alternative most
convenient. Place the command \Edges{...} in your style file (or
header) since its effect is is global. Other possible edge values
are 10 and 1000.

  If you use the command \Edges{...} at all, do so with care.  For
if you accidentally delete an \Edges{...} command your labels will
abruptly be badly misplaced and may logically but mysteriously
generate "dimension too big" errors under TeX and "off page" errors
under your driver.  

  You can dictate the edgescale for an individual figure by giving
the scale in brackets immediately after \SetLabels.  Thus, to
import into an article using say \Edge{100} a figure labelled using
another edgescale, say the original 1-by-1 default, you can use
\SetLabels[1]...\endSetLabels.

GETTING IT DOWN PAT

     Complicated labeling deserves the same respect as
complicated mathematics.  Do not expect it to come out perfect the
first time!  What is needed in either case is a mechanism to
repeatedly typeset troublesome pieces.

     One mechanism is always available.  One does complicated
labelling in a separate "test" file involving just the figure being
labelled;  a texpert will know how to \dump TeX's current state as
a temporary format that restarts rapidly at each retry.  Usually,
one then pastes the completed labelled figure back into the main
TeX file, but, of course, one can also \input it as an auxiliary
file.

     If you do not have a TeXpert at handy, here is a first
approximation to an efficient setup. By deletions reduce a copy
of your article to just a few lines before and after the figure.
Now label the figure, and finally, copy and paste the labelled
figure to the original article. Then copy the next figure to label
into this testbed and repeat. The TeXpert can improve the  speed
at which TeX starts up, by compiling a format specifically for
your article; just one caution: best NOT include in the format
ephemeral details of setup like \Set<mydriver>ArtSpecials (from
boxedeps.tex because this reads  figure dimensions which you may
change during your work session.

     An improved mechanism to repeatedly typeset troublesome
pieces is now available on the Macintosh; it is called LinoTeX;
see the same ftp sources.  It could be set up on many types
of computer.

     Before using labelfig.tex to attach labels to a graphics
object inserted using boxedeps.tex or BoxedArt.tex, make it a
firm rule to carefully adjust the bounding box using the trimming
commands of these packages, and also at least tentatively scale
and position the object. Beware of changing the grid inadvertently
after the labels have been positioned.  For example, correcting
the bounding box of a PostScript graphics object can foul up the
labels by changing the coordinate grid to which the labels are
attached. This is particularly true for the trimming  commands of
boxedeps.tex and BoxedArt.tex. However, as noted already, change
of scale is much less disruptive, and modest adjustments should be
well tolerated.

     Sometimes the labels protrude so far from the bounding box
of a figure that the figure has to be repositioned.  Best do this
by ad hoc spacing, say using \hglue and \vglue; altering the
bounding box would create a vicious circle.

     Remember that you are responsible for preventing labels
from overlapping. You are responsible for all label typography
including size and style. A label is really just about anything
that can be put in a TeX box. Note that spaces at the beginning
and end of labels will normally be suppressed; if you really want
them you must protect them with TeX braces.

     This package temporarily sets the \mathsurround parameter
of TeX to zero  while the labels are being affixed. This is done
because nonzero \mathsurround space would influence the position
of left and right aligned labels; then, when a texpert or printer
modifies mathsurround, diagram labeling might be disastrously
altered. There is a small price to pay involving labels that are
formatted as caption boxes including mathematics: you  may want or
need to specify an explicit mathsurround space within the caption
box; it will not influence anything outside.

     Those hostile to the use of * as separator between
the X and Y coordinates of label insertion points, are free to
impose another using \SetXYSeparator{<the new separator>}.  
Americans may prefer "," to "*" since they never use a 
comma as a decimal point; on the other hand, * may be more visible.

APPENDIX (I)  MERGING labelfig.tex LABELS INTO AN EPSF GRAPHICS OBJECT.

     As promised in the introduction, here is a recipe useful for
publishers. It works at least on Macintosh and at least for vectorized
graphics and Adobe type1 fonts.  (There is surely a similar recipe for
PCs under MSWindows.)

 (a)  Use boxedeps.tex utility to integrate the figure given by the eps
file, "x.eps" say, with a visible frame around it.  See
\ShowDisplacementBoxes command in boxedeps.tex.  To get precise results
automatically it is important to use the \Trim... commands of
boxedeps.tex making the "DisplacementBox" neatly fit the figure.

 (b)  Use the TeX printer driver and LaserWriter (versions >= 8.1.1) to
export to an EPSF the DVI page containing the integrated, labelled
figure. You now have an EPS file  "xx.eps"  that contains too much, and at
the wrong scale, and at wrong position.

 (c)  Convert the EPSF to an Adode Illustrator format EPSF using
the shareware utility called epsConvert by Sam Weiss
1993-- (currently $25).

 (d)  In Illustrator (or a compatible program), group the labels and the
"DisplacementBox"; copy them to the clipboard and paste them into "x.ps".
This step requires that all the label fonts be "visible to the Macintosh.

 (e)  Translate and scale the pasted group consisting of the labels plus
the "DisplacementBox" so as to make the "DisplacementBox" the bounding
box of (labelless) figure represented by "x.eps".  At this point the
labels will be correctly placed on the figure "x.eps".

 (f)  Ungroup and delete the "DisplacementBox".  The result is the
desired single EPS file, "x+.eps" say, It contains the original figure
plus its labels.  

     Using grouping and ungrouping appropriately in "x+.eps", a
publisher's staff can very efficiently improve label positions etc.

APPENDIX II)  SOME EXOTIC APPLICATIONS

     The grid of labelfig.tex is analogous to a light-table in
classical page makeup with wax or latex glue.  In principle, you
can use it to compose any page from its indivisible parts.  This
even has some of the artisanal charm of classical paste-up
provided you have a fast screen preview to make the process
"interactive".

     In practice labelfig.tex is a tool for nonstandard jobs.
Here are a few going beyond the labelling already discussed.

(I)  GRAPHICS INTEGRATION.

     This is accomplished by treating the imported graphics
objects as labels.  The underlying graphics object is then
typically an empty  \vbox to <dimension>{\vfill} in a TeX
\midinsert...\endinsert construction.  A label line
might be of the form

   (.1*.1) \\

The exact form of the special command varies from driver to
driver.  However, in the case of encapsulated PostScript graphics
(EPSF norm), by relying on boxedeps.tex, one can have the
following standard syntax (independant of driver  (see
boxedeps.doc for details.
  
  (.1*.1) \BoxedEPSF{MyFigure scaled <scale in mils>}\\

This may be slow since it requires TeX to read the PostScript
file to read bounding box using many complex macros.  So you
may want to try

  (.1*.1) \EPSFSpecial{MyFigure}{<scale in mils>}\\

which is fast and driver independant, but it squashes the
bounding box, normally to its lower left corner.

     Similarly for graphics of the Macintosh PICT norm ---
using BoxedArt.tex (same sources) in place of boxedeps.tex.

     This approach to integration is to be recommended when
one is assembling a composite graphics object.

 (II)  COMMUTATIVE DIAGRAM ENHANCEMENT

     Commutative diagrams or arrays of mathematical objects
connected by arrows of various sorts are common in mathematics.
The mathematical objects require the use of TeX.  Recently TeX
acquired a good collection of arrows of all slopes --- that of
LamSTeX --- plus pwerful macros to build the diagrams.

     However, even the LamSTeX collection is often
inadequate; it lacks for example double shafted arrows, dotted
arrows and curved arrows. Fortunately it is possible to produce
such arrows on an individual basis using sophisticated graphics
programs such as Illustrator and AldusFreehand (both serving
the EPSF norm) or using Metafont (with its public domain norm).
Since the creation of each new arrow is a work of love, you
probably want to limit the number of arrows by using LamSTeX
for most arrows. The 40K commutative diagram module of LamSTeX
has been adapted to work with AmSTeX and a copy may be posted
with LabelFig and related files. Unfortunately no one has yet
offered a version that works with Plain TeX or LaTeX.

       Suffice it here to say that when the exotic arrow has
been somehow imported into TeX, labelfig.tex treats it as a
label that one affixes to the commutative diagram.  Two other
steps will be treated in separate notes, namely the matter of
extracting the dimension specifications for the arrow and the
construction of the arrow --- for these steps are far from
unique and often depend intimately on your computer environment. 
Notes for the Macintosh-Textures-Illustrator combination are
found in the file ExoticArrows.doc.

 (III) NESTING 

Ingenuity pays off in exploiting labelfig.tex. One can
mix graphics and typography quite freely.  labelfig.tex is good
for freeform or overlapping arrangements, while boxedeps.tex (or
BoxedArt.tex) is best for regimented non-overlapping
arrangements --- and the two can be combined.

     The default behavior of labelfig.tex is not ideal 
for nesting objects, because to prevent trouble for beginners
the register for labels is globally cleared when \AffixLabels
concludes.  But there are switches available

      \LabelsGlobal      \LabelsLocal

which change this.  To understand this, extend the above test 
by something like:

 \LabelsLocal
 \SetLabels
    (.5*.5) AAA\\
 \endSetLabels

 {%%% Watch for influence of braces!!
 \SetLabels
    (.5*.5) ZZZ\\
 \endSetLabels
   \AffixLabels{\FirstQuadrant}
 }

   \AffixLabels{\FirstQuadrant}

     There are however potential pitfalls.  Neither
labelfig.tex nor boxedeps.tex has been tested under extreme
conditions. Problems may occur if their procedures are
indiscriminately nested. For boxedeps.tex (not labelfig.tex)
there is a precise cause for worry, namely many of its
variables are "global", which means that TeX braces will not
provide the protection one might expect.

COMMAND SUMMARY FOR labelfig.tex

  Here [...] means optional (one or zero)
       [...]* means any number of such constructs

  \SetLabels
    [[<P>](<X><Sep><Y>) <label> \\]*
  \endSetLabels
  \ShowGrid  % this makes the grid visible (once)
  \AffixLabels{<the figure>}

   --- <P> is tack position, one of eleven or empty
              order irrelevant

                   \L\T      \T      \R\T

                   \L\E      \E      \R\E

                     \L               \R

                   \L\B      \B      \R\B

   --- (<X><Sep><Y>) insertion point;
  <Sep> is separator, = * by default;
  \SetXYSeparator{<Sep>} changes it.
   <X> and <Y> are real numbers

  --- <label> a label to attach 

  --- <the figure> the figure to label 

  \GlobalLabels (default)     
  \LocalLabels  setting for nested constructs.

 \Grids makes ALL grids appear; \HideGrid then makes just next disappear.
 \noGrids returns to default.  The commands are always global.

 \GridLineWidth{<dimension>} adjusts width of grid lines. Default is very
small, to give "hairline" effect. If your grid lines are missing try
setting \GridLineWidth{1pt}.

 \Edges#1 globally changes the edge size of all grids to the numerical 
value #1, which must be 1, 10, 100, or 1000.  The default is 1.

VERSION HISTORY.
 --- Jan 1993: \Edges#1 and [??] option after \SetLabels
 --- July 1992: \Grids, \noGrids, \HideGrid;
       Gridlines become hairlines; \GridLineWidth{<dimension>}.
 --- Oct 1991, Jan 1992: \SetXYSeparator{<Sep>},  \LabelsGlobal,
       \LabelsLocal.
 --- July 1991: first release

Address for bugs and other feedback:

        Raymond S\'eroul
        IREM and Lab. de Typographie Informatise
        Univ. Rene Descartes
        Strasbourg

    Tel 33-88-41-63-45
    Email:  A18645@FRCCSC21.BITNET

        Laurent Siebenmann
        Mathematique, Bat. 425,
        Univ de Paris-Sud,
        91405-Orsay,
        France

    Tel 33-1-6941-7949; 
    Email: lcs@topo.math.u-psud.fr

\newcommand\C{{\mathbb C}}

\newcommand\D{{\mathbb D}}

\newcommand\Z{{\mathbb Z}}
\newcommand\R{{\mathbb R}}

\newcommand\co{\colon}
\renewcommand\Im{\mathop{\mathrm{Im}}}
\renewcommand\Re{\mathop{\mathrm{Re}}}

\newcommand\OC{{\overline{\mathbb C}}}
\newcommand\OR{{\overline{\mathbb R}}}
\newcommand{\CC}{\OC}
\DeclareMathOperator{\arsinh}{arsinh}

\newtheorem{corollary}{Corollary}

\newtheorem{theorem}{Theorem}

\newtheorem{lemma}{Lemma}

\begin{document}
\author{Lukas Geyer}
\address{Lukas Geyer\\University of Michigan\\Department of Mathematics\\
  525 E.~University\\Ann Arbor, MI 48109\\USA}
\email{lgeyer@umich.edu}
\thanks{The first author was supported by the Alexander von Humboldt
  Foundation. The second author was supported by NSF grants DMS-0400636, 
 DMS-0244421, and DMS-0244547.}
\author{Sergei Merenkov}
\address{Sergei Merenkov\\ University of Michigan\\Department of Mathematics\\
  525 E.~University\\Ann Arbor, MI 48109\\USA}
\email{merenkov@umich.edu}
%\thanks{The second author was supported by NSF grants DMS-0400636, 
%DMS-0244421, and DMS-0244547.}
%\thanks{}
%\address{Department of Mathematics\\ Weizmann Institute of Science\\
%Rehovot 76100, Israel} \email{itai@math.weizmann.ac.il}
%\address{Department of Mathematics\\ Purdue University\\ West Lafayette,
%Indiana 47907} \email{smerenko@math.purdue.edu}
%\address{Microsoft Research\\ One Microsoft Way\\ Redmond, WA 98052}
%\email{schramm@microsoft.com}

\abstract{We prove the existence of a hyperbolic surface 
spread over the sphere for which the projection map has all
its singular values on the extended real line, and such that the preimage
of the extended real line under the projection map is homeomorphic
to the square grid in the plane. This answers a question raised by 
{\`E}.~B.~Vinberg.}
\endabstract

\maketitle

\section{Introduction}\label{S:Intro}

According to the uniformization theorem, an open simply connected
Riemann surface is conformally equivalent to either the complex plane
$\C$, or the open unit disc $\D=\{z\in\C:|z|<1\}$. In the former case
it is said to be \emph{parabolic}, and in the latter it is called
\emph{hyperbolic}.  We are interested in the application of the
uniformization theorem to clasically defined surfaces spread over the sphere.

A {\emph{surface spread over the sphere}} is a pair $(X, \psi)$, where
$X$ is a topological surface and $\psi\co X\to\OC$ a continuous, open and
discrete map.  Here $\OC$ denotes the Riemann sphere.  The map $\psi$ is
called a {\emph{projection}}. Two such surfaces $(X_1, \psi_1),\ (X_2,
\psi_2)$ are {\emph{equivalent}}, if there exists a homeomorphism $\phi\co
X_1\to X_2$, such that $\psi_1= \psi_2\circ\phi$.  According to a theorem of
Sto\"{\i}low \cite{sS56}, there exists a unique conformal structure on
$X$ such that $\psi$ is a holomorphic map. We call this the
\emph{pull-back structure} on $X$.  If $X$ is open and simply connected, 
then as a Riemann surface with the pull-back structure it is conformally 
equivalent to either $\C$ or $\D$.
Equivalent surfaces have the same conformal type.

A continuous, open and discrete map $\psi$ 
near each point $z_0$ is homeomorphically
equivalent to a map $z\mapsto z^k$. The number $k=k(z_0)$ is called
the {\emph{local degree}} of $\psi$ at $z_0$. If $k\neq 1$, $z_0$ is
called a {\emph{critical point}} and $\psi(z_0)$ a {\emph{critical
value}}. The set of critical points is a discrete subset of $X$.
A point $a\in\OC$ 
is called an \emph{asymptotic value} of $\psi$, if there exists 
a curve $\gamma\colon [0, t_0)\to X$ such that
$$
\gamma(t)\to\infty {\text{ and }} 
\psi(\gamma(t))\to a {\text{ as }} t\to t_0.
$$ 
A point $a\in \OC$ is called a \emph{singular value} of
$\psi$ if $a$ is either a critical or an asymptotic value of $\psi$. 

{\`E}.~B.~Vinberg~\cite{eV89} introduced the interesting class
$\mathcal{V}$ of surfaces spread over the sphere $(X, \psi)$, such that
all singular values of $\psi$ are contained in $\OR=\R\cup\{\infty\}$.  For such
mappings $\psi$, the components of the preimages of the upper and lower
hemispheres are called \emph{cells}.  The preimage of
the extended real line is a disjoint union of the critical points,
called \emph{vertices}, and open arcs, called \emph{edges}.
Each vertex has 
an even number (which is greater than or equal to four)
of edges emanating from it. The map $\psi$ in this case is called a 
\emph{cellular map}. 
The set
$\psi^{-1}(\OR)$, understood as an embedded graph $N$ in
$X$, determines the topological properties of $\psi$. Such an embedded
graph $N$ is called a \emph{net}. Suppose that the edges
of this graph are labeled by positive numbers such that the sum over
the edges of each cell is $2\pi$.  If $\OC$ is endowed with the
spherical metric, we assume that $\psi$ maps an edge labeled by
$\alpha$ onto an arc of length $\alpha$.  Since the length of the real
line in the spherical metric is $2\pi$, this is possible.
Equivalently, one can label the vertices of a net and 
ends of edges going to infinity by points of $\OC$,
so that the order of labels correspond to the natural order on
$\OR$.

It is an interesting question in which cases the conformal type of $X$
is determined by the combinatorics of a net, independent of the
labeling.  G.~MacLane~\cite{gM47} and Vinberg~\cite{eV89} 
considered the periodic net
$\cos^{-1}(\R)$.  They showed that all 
labelings of this net by real numbers and $\infty$, with ends of edges 
going to infinity labeled by $\infty$, produce parabolic surfaces.

In his paper, Vinberg asked whether the conformal type is determined
by combinatorics in the case of the square grid net. It is clear that if
the labels on the vertices of every two adjacent square cells are
symmetric with respect to their common edge, the resulting surface is
parabolic, and the corresponding cellular map is a fractional linear
transformation of some Weierstrass elliptic function. The main result
of this paper is the following theorem, which answers the question of 
Vinberg. 
\begin{theorem}\label{T:Maintheorem}
  There exists a hyperbolic surface of class $\mathcal V$ whose net is
  homeomorphic to the square grid.
\end{theorem}

{}From the point of view of the type problem, the class $\mathcal V$ is
a generalization of a classically studied class $\mathcal S$.  A
surface spread over the sphere $(X, \psi)$ is said to belong to the
class $F(a_1, \dots, a_q)$, if $\psi$ restricted to the complement
of $\psi^{-1}\bigl(\{a_1,\dots,a_q\}\bigr)$ is a covering map onto its
image $\OC\setminus\{a_1,\dots,a_q\}$.  One defines $F_q :=
\displaystyle \bigcup_{a_1,\ldots,a_q} F(a_1,\ldots,a_q)$ and $\mathcal S := \bigcup_q F_q$.

Surfaces of class $\mathcal S$ have a combinatorial representation
in terms of so-called labeled Speiser graphs.
Assume that $(X,\psi)\in F_q$ and $X$ is open.
We fix a Jordan curve $L$, visiting the points $a_1,\dots, a_q$
in cyclic order.
The curve $L$ is usually called a {\emph{base curve}}. It
decomposes the sphere into two simply connected regions $H_1$,
the region to the left of $L$, and $H_2$, the region to the right of
$L$.
Let $L_i$, $i=1,2,\dots,q$, be the arc of $L$
from $a_i$ to $a_{i+1}$ (with indices taken
modulo $q$). Let us fix
points $p_1$ in $H_1$ and $p_2$ in $H_2$, and choose $q$
Jordan arcs $\gamma_1,\dots, \gamma_q$ in ${\OC}$, such that each
arc $\gamma_i$ has $p_1$ and $p_2$ as its endpoints, and has a
unique point of intersection with $L$, which is in $L_i$.
We take these arcs to be interiorwise disjoint, that is,
$\gamma_i\cap\gamma_j=\{p_1,p_2\}$ when $i\neq j$.
Let $\Gamma'$ denote the graph embedded in $\OC$, whose vertices are
$p_1,\ p_2$, and whose edges are $\gamma_i,\ i=1,\dots, q$, and let
$\Gamma=\psi^{-1}(\Gamma')$. We identify $\Gamma$ with its image in
$\R^2$ under an orientation preserving homeomorphism of $X$ onto $\R^2$.
%Clearly it does not depend on the choice of the points $p_1,\
%p_2$, and the curves $\gamma_i,\ i=1, \dots, q$.
The graph $\Gamma$ has the following properties:
It is infinite, connected, homogeneous of degree $q$, and bipartite.
A graph, properly embedded in the plane and having these properties
is called a {\emph{Speiser graph}}, also known as a
{\emph{line complex}}. The vertices of a Speiser graph $\Gamma$
are traditionally marked by $\times$ and $\circ$, such that each edge
of $\Gamma$ connects a vertex marked $\times$ with a vertex marked $\circ$.
Each {\emph{face}} of $\Gamma$, i.e. a connected component of
$\R^2\setminus\Gamma$, has either a finite even number of edges along its
boundary, in which case it is called an
{\emph{algebraic elementary region}},
or infinitely many edges, in which case it is called a
{\emph{logarithmic elementary
region}}. Two Speiser graphs $\Gamma_1, \ \Gamma_2$ are said to be
{\emph{equivalent}}, if there is a sense-preserving homeomorphism
of the plane, which takes $\Gamma_1$ to $\Gamma_2$.
% Below we refer to an equivalence class as a Speiser graph.

The above construction is reversible. Suppose that the faces of a
Speiser graph $\Gamma$ are labeled by $a_1, \dots, a_q$, so that when
going counterclockwise around a vertex $\times$, the indices are
encountered in their cyclic order, and around $\circ$ in the reversed
cyclic order.  We fix a simple closed curve $L\subset\OC$ passing
through $a_1, \dots, a_q$.  Let $H_1,H_2,L_1,\dots,L_{q}$ be as
before.  Let $\Gamma^*$ be the planar dual of $\Gamma$.  If $e$ is an
edge of $\Gamma^*$ from a face of $\Gamma$ marked $a_j$ to a face of
$\Gamma$ marked $a_{j+1}$, let $\psi$ map $e$ homemorphically onto the
corresponding arc $L_j$ of $L$.  This defines $\psi$ on the edges and
vertices of $\Gamma^*$.  We then extend $\psi$ to map the faces of
$\Gamma^*$ homeomorphically to $H_1$ and $H_2$,
respectively. This defines a surface spread over the sphere $(\R^2,
\psi)\in F(a_1, \dots, a_q)$. The conformal structure is induced
by the conformal structure of $\OC$ if one considers a copy
of $H_1$ or $H_2$ for each face of $\Gamma^*$, and glues $H_1$ to $H_2$
using the identity map whenever the corresponding faces are adjacent. 
See~\cite{Nevbook} for further details.

The proof of Theorem~\ref{T:Maintheorem} is in two parts. In 
Section~\ref{S:HypS} we construct a hyperbolic surface of class $\mathcal S$
whose net is a certain degeneration of the square grid, namely some
of the sides of the grid are collapsed to a point. In Section~\ref{S:unrav}
we use quasiconformal deformations to obtain a surface 
of class $\mathcal V$ with the desired net. For a background on the theory
of quasiconformal maps see~\cite{LV}.

We would like to mention two interesting open questions connected to
this result. First, is it possible to find a hyperbolic surface spread over
the sphere with combinatorics of the square grid, which has a labeling 
that is symmetric with respect to a line? 
This would correspond to a real meromorphic map in the unit
disc such that the preimage of the extended real line is homeomorphic
to the square grid. The second question is whether MacLane-Vinberg's result
still holds if one allows non-symmetric labeling of the cosine net,
i.e.~whether the cosine net always gives a parabolic surface 
in the case when we
do not assume that labels on the ends of edges going to infinity
in the upper half-plane are the same as in the lower half-plane.

{\sc Acknowledgments.} We would like to thank Alexandre Eremenko for
suggesting this problem and Mario Bonk for helpful discussions. We
would also like to thank the referee for valuable comments and for
suggesting a simplification of the proof in Section~\ref{S:HypS}.

\section{A hyperbolic surface of class $\mathcal S$}\label{S:HypS}

\begin{figure}
\centerline{\AffixLabels{%
\includegraphics*[height=1.9in]{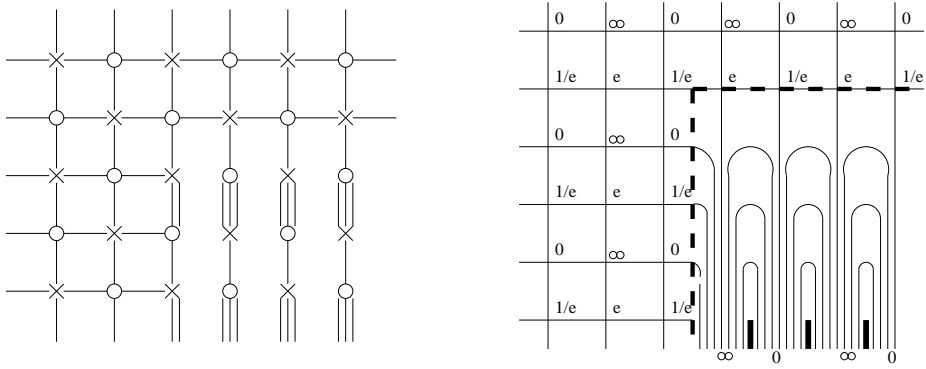}
}}
\caption{\label{f.Speiser}Speiser graph $\Gamma$ and the corresponding net.}
\end{figure}

Consider a planar graph $G$ whose vertices form the set 
$\{(m,n)\co\ m,n\in\Z\}$, and whose edge set is 
\begin{align}
&\{[(m,n),(m,n+1)]\co\ m,n\in \Z\}\cup \{[(m,n),(m+1,n)]\co\ m,n\in \Z, 
n\geq0\}\notag\\
&\cup \{[(m,n),(m+1,n)]\co\ m,n\in \Z, 
n<0, m<0\}.\notag
\end{align}
We denote by $\circ$ every vertex $(m,n)$ with $|m+n|$ even,
and by $\times$ every vertex $(m,n)$ with $|m+n|$ odd. By replacing edges 
of the set 
$$\{[(m,-2n), (m, -2n+1)]\co\ m,n\in\Z, m\geq0, n\geq1\}
$$
in the graph 
$G$ by multiple edges of the same form, we obtain a homogeneous graph
of degree four. This is a Speiser graph which we denote by $\Gamma$
(see Figure~\ref{f.Speiser}).
We label its faces by $0, 1/e, e, \infty$, so that this order is preserved
when going counterclockwise around $\times$, and the face adjacent to
the vertex $(0,0)$ which is in the first quadrant is labeled by $e$.

The Speiser graph $\Gamma$ with so prescribed labeling defines a unique
surface $(X, \psi)$ of class $\mathcal S$, or, more precisely, of class
$F_4(0,1/e,e,\infty)$. Here we assume that the base curve is the extended 
real line. The rest of this section is devoted to the proof 
that the surface $(X,\psi)$ is of hyperbolic type.

We will show hyperbolicity by first cutting the surface into two parts $A$
and $B$, and finding explicit quasiconformal maps to horizontal strips. We
then calculate the asymptotics of the gluing maps and use a
Theorem of L.~I.~Volkovyski\u{\i}~\cite{lV50} 
to show hyperbolicity of the resulting surface.

Let $B$ be the subsurface of $X$ with boundary which is obtained by
gluing together hemispheres corresponding to the vertices of the
Speiser graph $\Gamma$ from the set $\{ (m,n) \co\ m>0, \, n\leq 0 \}$, 
together with quarterspheres from the set $\{ (0,n) \co\ n\leq 0\}$. 
Quarterspheres here are the parts of the hemispheres outside the
unit disc. Let $A$ be the complement of $B$ in $X$. Then the interiors
of $A$ and $B$ are simply connected and the common boundary is a
simple curve (see Figure~\ref{f.Speiser}).

%\subsection{Uniformization of $A$ and $B$}

\begin{figure}
\centerline{\AffixLabels{%
\includegraphics*[height=2.1in]{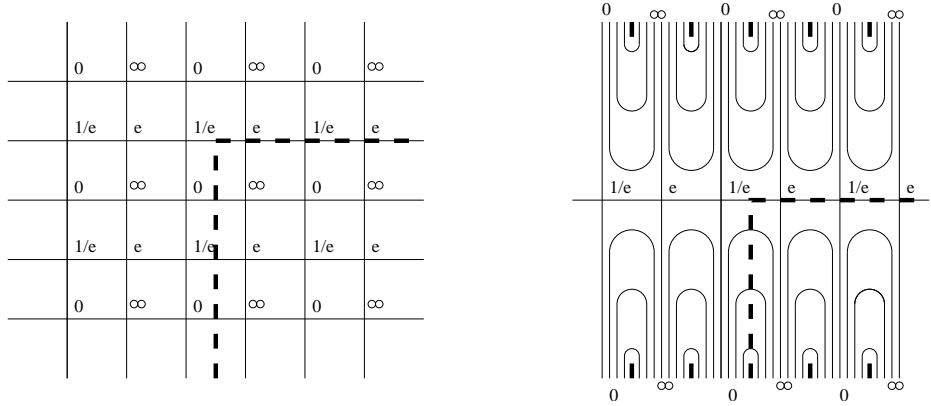}
}}
\caption{\label{f.Uniformized}Uniformizations of $A$ and $B$ .}
\end{figure}

The part $A$ is uniformized by the strip $S_1=\{z\co 0<\Im z<\pi\}$,
using the map
\begin{equation*}\label{E:B}
  \wp (\alpha(\exp(z))),
\end{equation*}
where $\alpha$ is the quasiconformal map given in polar coordinates by
$(r,\theta) \mapsto (r,3\theta/2)$, $0<\theta < \pi$, and $\wp$ is the conformal map of a rectangle
$\{ x+iy: 0<x<\pi/2, \, 0<y<\tau \}$ to the domain $\{ |z|>1, \, \Im z >0 \}$,
with $\wp(0)=1$, $\wp(\pi/2)=e$, $\wp(\pi/2 + i\tau) = \infty$ and $\wp(i\tau) = -1$. By
reflection $\wp$ extends to an elliptic function in the plane with
periods $2\pi$ and $2i\tau$. The map $\alpha \circ \exp$ maps $S_1$ to 
quadrants~I--III, and $\wp$ maps these quadrants to $A$ (see
Figure~\ref{f.Uniformized}).
The map $\wp$ sends the real line to the real
line, and the imaginary axis to the unit circle. The critical points
are at $\pm\pi/2$, $\pm\pi/2+i\tau$ and all equivalent points. The critical
values are $\wp(\pi/2)=e$, $\wp(-\pi/2)=1/e$, $\wp(\pi/2+i\tau)=\infty$ and
$\wp(-\pi/2+i\tau)=0$.

The part $B$ is uniformized by the strip $S_2=\{z\co \pi< \Im z < 2\pi \}$,
using the map 
\begin{equation*} \label{E:A}
\exp(\sin(\beta(\exp(z)))),
\end{equation*}
where $\beta$ is the quasiconformal map given in polar coordinates by
$(r,\theta) \mapsto (r,\theta/2)$, $-\pi < \theta < 0$. The map $\beta 
\circ \exp$ maps $S_2$ to quadrant~IV, and $\exp \circ \sin$ maps this 
quadrant to $B$ (see Figure~\ref{f.Uniformized}).

The gluing map $f$ from the boundary $\{\Im z=0\}$ of $S_1$ to the boundary
$\{\Im z=2\pi \}$ of $S_2$ and
the gluing map $g$ from the boundary $\{\Im z=\pi \}$ of $S_1$ to the boundary
$\{\Im z=\pi \}$ of $S_2$ are given by
$$
f(t)=\log(e^t+p(e^t)),\ \ g(t)=\log\bigg(-\arsinh\bigg(\frac{\pi}{\tau}e^t+
q(e^t)\bigg)\bigg),
$$
where $p(x)$ and $q(x)$ are real-analytic periodic 
functions with derivatives bounded 
below by $c>-1$. Differentiating, we find that 
$f'(t)\geq(c+1)/2$ and $g'(t)\leq C/t$ for $C>0$ and $t$ large enough. 
By \cite[Theorem 24, p.~89]{lV50}  (see also the remark following it),
the resulting surface is hyperbolic.

\section{Unraveling logarithmic staircases}\label{S:unrav}
In this section we will use quasiconformal deformations of $(X,\psi)$ to
obtain a surface of the same type with the net homeomorphic to the
square grid.

%Let $(X,\psi)$ be a Riemann surface spread over the sphere. 
A simply
connected subsurface $(Y,\psi)$ with $Y\subset X$ is a \emph{logarithmic
  staircase} over $a\in \R$ in $(X,\psi)$ if there exists $c\in\R$ such that
$\psi|_Y : Y\to D_{e^c}^*(a)$ is a regular covering map, i.e. there exists
a conformal map $\phi:Y\to\{\Re z < c\}$ such that $\psi(z) = a+ e^{\phi(z)}$ for $z\in
Y$.  The preimage of the real line consists of infinitely many
half-lines $l_k = (-\infty,c) + i\pi k$, $k\in\Z$, in the $\phi$-coordinate. A
simply connected subsurface $(Y,\psi)$ is a logarithmic staircase over
$\infty$ if $(Y,1/ \psi)$ is a logarithmic staircase over 0.
The surface $(X,\psi)$ contains infinitely many logarithmic staircases
over $0$ and $\infty$.  

\begin{figure}
\centerline{\AffixLabels{%
\includegraphics*[height=1.9in]{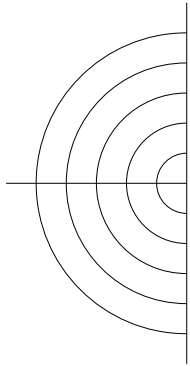}
}}
\caption{\label{f.approx}Topological picture of the net of the
  $cos$-spine.}
\end{figure}

In order to obtain a surface whose net is
homeomorphic to the square grid, we will replace all of those
logarithmic staircases by \emph{cosine spines}.
We say that $(Y,\psi)$ is a \emph{$\cos$-spine} in $(X,\psi)$ if $Y\subset X$ is simply connected, $\psi$ is a cellular map, 
and $(\psi|_Y)^{-1} (\R)$ is a union of infinitely many
lines $l_0$, $l_1$,\ldots, where $l_0$ has one endpoint on $\partial Y$, all other
lines have two endpoints on $\partial Y$ and $l_0$ intersects every $l_k$ once,
and none of the other $l_k$ intersect (see Figure~\ref{f.approx}). The
line $l_0$ is the \emph{axis of symmetry} of the $\cos$-spine.

\begin{theorem}
  \label{thm:app}
  Let $(X,\psi)$ be a Riemann surface spread over the sphere and $(Y,\psi)$
  a logarithmic staircase over some $a\in\OR$ in $(X,\psi)$. Let $z_0$ be a
  point on $\partial Y$ with $\psi(z_0)\in \R$.  Endow $X$ with the pull-back
  complex structure. Then there exists a universal constant $K>1$ 
  and a $K$-quasiregular map
  $\tilde{\psi} : X \to \CC$ with $\tilde{\psi} = \psi$ outside $Y$ and
  $(Y,\tilde{\psi})$ being a $\cos$-spine whose axis of symmetry passes
  through $z_0$.
\end{theorem}

An immediate corollary is the following.
\begin{corollary}\label{C:spines}
  \label{cor:app}
  Let $(X,\psi)$ be a Riemann surface spread over the sphere and
  $(Y_1,\psi)$, $(Y_2,\psi)$,\ldots be mutually disjoint logarithmic staircases
  in $(X,\psi)$. Furthermore, let $z_j$ be a point on $\partial Y_j$ with
  $\psi(z_j) \in \R$.  Then there exists a Riemann surface $(X,\tilde{\psi})$
  of the same type as $(X,\psi)$, where $\tilde{\psi}=\psi$ outside $\bigcup_j Y_j$,
  and $\tilde{\psi}|_{Y_j}$ is a $\cos$-spine for each $j$ with its axis
  of symmetry passing through $z_j$.
\end{corollary}

{\sc Proof.} If we apply Theorem \ref{thm:app} to each logarithmic
staircase separately, we obtain a quasiregular map $\tilde{\psi}:X \to \CC$
with the desired net. This implies that the conformal structures
induced by $\psi$ and $\tilde{\psi}$ are quasiconformally equivalent. The
type of a simply connected Riemann surface is preserved under
quasiconformal homeomorphisms. \qed

Before proving Theorem \ref{thm:app} we need some preliminary results.

Let $S$ denote the strip $\{0<\Im z<1\}$.

\begin{lemma}
  \label{lemma:interpolation}
  Let $f_0,f_1:\R\to\R$ be diffeomorphisms and $M\geq 1$ with $1/M \leq
  f_k'(x) \leq M$ and $|f_0(x)-f_1(x)|\leq M-1$ for $k=0,1$ and all $x\in\R$.
  Then there exists a $K(M)$-quasiconformal map $f:S\to S$ with $f(x+ik)
  = f_k(x)+ik$ for $k=0,1$. Furthermore, $\lim_{M\to 1} K(M) = 1$. 
\end{lemma}
{\sc Proof.} Define $f(x+iy) := (1-y) f_0(x) + y f_1(x) + iy$ in
$S$. The fact that $f_0$ and $f_1$ are increasing homeomorphisms of
the real line immediately implies that $f$ is a homeomorphism of $S$
onto itself, preserving boundary components.
The partial derivatives are
\begin{eqnarray*}
  f_x(x+iy) &=& (1-y) f_0'(x) + y f_1'(x), \\
  f_y(x+iy) &=& f_1(x) - f_0(x) + i.
\end{eqnarray*}
Viewing the complex plane as $\R^2$, the matrix of the differential is
\[
Df(x,y) = 
\begin{pmatrix}
  (1-y) f_0'(x) + y f_1'(x) & f_1(x) - f_0(x). \\
  0 & 1
\end{pmatrix}
\]
Thus $Jf(x,y) = (1-y) f_0'(x) + y f_1'(x) \geq 1/M$ and $\| Df(x,y) \| \leq
\max (1,(1-y)f_0'(x) + y f_1'(x) + |f_1(x)-f_0(x)|) \leq 2M-1$. For the
dilatation these estimates give
\[
Kf(x,y) = \frac{\|Df(x,y)\|^2}{Jf(x,y)} \leq M(2M-1)^2,
\]
which is finite for all $M$ and goes to 1 for $M\to 1$, as claimed. \qed

\begin{figure}
\centerline{\AffixLabels{%
\includegraphics*[height=1.9in]{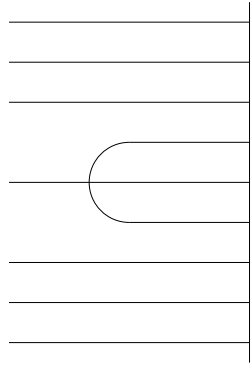}
}}
\caption{\label{f.Onecrit}Topological picture of the net of the
  function in Lemma~\ref{lemma:onecrit}.}
\end{figure}

Let $D\subset \C$ be a domain which is symmetric w.r.t.~the real line.
A function $G:D \to \C$ is called \emph{real} if $G(\overline{w}) =
\overline{G(w)}$ for all $w\in D$.

\begin{lemma}
\label{lemma:onecrit}
There exists a real quasiregular function $G: \{ \Re w \leq 1 \} \to \C$ and a
constant $M>0$ such that
  \begin{enumerate}
  \item $G(w) = e^w$ for $\Re w =1$ and for $\Re w \leq -M-1$
  \item $G$ has one real critical point and the preimage of the real
    line consists of infinitely many mutually disjoint rays to $\infty$
    that start on the line $\{ \Re w = 1 \}$, and an arc which has two
    endpoints on that line, intersects the real axis once, and is
    disjoint from all other rays (see Figure \ref{f.Onecrit}).
  \end{enumerate}
\end{lemma}

{\sc Proof.} 
Define
\[
H(w) := \frac{w+4}{w-4} e^w.
\]
We want to interpolate between $H$ and $\exp$ (locally)
quasiconformally in the vertical strips $S'= \{0 < \Re w < 1 \}$ and
$S'' = \{ -M-1 < \Re w < -M \}$, where $M$ is to be chosen. For
interpolation in $S'$ we will use Lemma \ref{lemma:interpolation} and
thus we have to check the conditions for the lifts (logarithms) of the
boundary mappings of $H$ on $\{ \Re w = 0 \}$ and $\exp$ on $\{ \Re w = 1 \}$.
In our case $f_1$ is the identity mapping as the imaginary part of the
logarithm of $\exp(1+ix)$. For $f_0$ we get
\begin{eqnarray*}
  f_0(x) &=& \Im \log H (ix) \\
  &=& \Im \left(\log \frac{4+ix}{-4+ix} + ix\right) \\
  &=& x + \arctan \frac{8x}{16+x^2}.
\end{eqnarray*}
{}From the boundedness
of the arctangent we immediately have $|f_1(x)-f_0(x)| \leq \frac\pi2$. 
Furthermore,
\begin{eqnarray*}
  f_0'(x) 
% &=&
% 1 + \frac{1}{1+\left(\frac{8x}{16+x^2}\right)^2}
%    \frac{8(16+x^2) - 8x \cdot 2x}{(16+x^2)^2} \\
    &=& 1+ \frac{8(16 - x^2)}{(16+x^2)^2+64x^2},
\end{eqnarray*}
from which we get 
\[
  |f_0'(x) - 1| = \frac{|8(16 - x^2)|}{|(16+x^2)^2+64x^2|} \leq  \frac12
\]
for all $x\in\R$. An application of Lemma 1 shows that there exists a
quasiconformal map $\tilde{G}$ in $S'$ with boundary values $\log
H(w)$ on $\{ \Re w = 0\}$ and $w$ on $\{ \Re w = 1 \}$. Then $G = \exp \circ
\tilde{G}$ is quasiregular in $S'$ with boundary values $H$ and
$\exp$, respectively.

For interpolation in $S''$,
% we can not use Lemma
%\ref{lemma:interpolation} because $H$ does not map the line $\{ \Re w =
%-M \}$ to a circle. 
fix an increasing $C^\infty$-function $\eta:\R \to [0,1]$
with $\eta(t) = 0$ for $t<-1$ and $\eta(t)=1$ for $t>0$. Define
\begin{eqnarray*}
 G(w) &=& \eta(M+ \Re w) H(w) + (1-\eta(M+\Re w)) e^w \\
 &=& e^w \left[ \eta(M+\Re w) \frac{w+4}{w-4} + 1-\eta(M+\Re w) \right] \\
 &=& e^w \left[ \eta(M+\Re w) \frac{8}{w-4} + 1 \right]
\end{eqnarray*}
for $\Re w \leq 0$. This agrees with the earlier definition of $G$ on the
imaginary axis. The function $G$ is smooth in the left half-plane, and
it is holomorphic outside $S''$. We claim that $G$ is locally
quasiconformal there when $M$ is sufficiently large. The partial
derivatives of $G$ in $S''$ are (we omit the argument $M+\Re w$ for $\eta$
and $\eta'$ here)
\begin{eqnarray*}
  \frac{\partial G}{\partial\overline{w}} &=& e^w  \eta' \frac{4}{w-4} 
\end{eqnarray*}
and
\begin{eqnarray*}
  \frac{\partial G}{\partial w} &=& e^w \left[ \eta \frac{8}{w-4} + 1 +
    \eta' \frac{4}{w-4} - \eta \frac{8}{(w-4)^2} \right],
\end{eqnarray*}
thus we have the estimates (assuming $M>4$)
\begin{eqnarray*}
  \left| e^{-w} \frac{\partial G}{\partial\overline{w}} \right| &\leq& \frac{4\| \eta' \|}{M-4},
\end{eqnarray*}
and
\begin{eqnarray*}
  \left| e^{-w} \frac{\partial G}{\partial w} \right| &\geq& 1- \frac{8}{M-4} - \frac{4
    \| \eta' \|}{M-4} - \frac{8}{(M-4)^2},
\end{eqnarray*}
where $\| . \|$ stands for the supremum norm. Choosing $M$ large enough, we can
ensure $\left| e^{-w} \frac{\partial G}{\partial\overline{w}} \right| \leq 1/3$ and
$\left| e^{-w} \frac{\partial G}{\partial w} \right| \geq 2/3$, which gives $\left|
  \frac{\partial G}{\partial\overline{w}} \left/ \frac{\partial G}{\partial w} \right. \right| \leq 1/2$
for the complex dilatation of $G$. This shows that $G$ is locally
quasiconformal in $S''$, thus it is a
quasiregular map.  \hfill\qed(Lemma~\ref{lemma:onecrit})

{\sc Proof of Theorem \ref{thm:app}.} Assuming that $a\in\R$, 
let $\phi:Y\to\{\Re w < c\}$ be a
uniformization of $Y$ with $\phi(z_0)\in \R$ and
$\psi(z)= a+e^{\phi(z)}$ for $z\in Y$. 
To obtain the desired quasiregular map we 
modify the exponential map in the half-plane $\{ \Re w < c\}$. 
Define
\[
\begin{split}
g(w) &:=
  e^{c-k(M+2)}G(w-c+k(M+2))\\
     &\text{ for } -M-1 < \Re w-c + k(M+2) \leq 1,
\end{split}
\]
where $k=0,1,2,\ldots$ and $G$ is the function from
Lemma~\ref{lemma:onecrit}. 

The map $\tilde{\psi}$ is given by 
$$
\tilde{\psi}(z)=\begin{cases} 
\psi(z),\ z\notin Y,\\
a+ g(\phi(z)),\ z\in Y.
\end{cases}
$$

If $a=\infty$, the function $\tilde{\psi}$ is defined in $Y$ by the formula
$\tilde{\psi}(z)=1/g(\phi(z))$.

\hfill\qed(Theorem~\ref{thm:app})

{\sc Proof of Theorem~\ref{T:Maintheorem} continued.} 
Theorem~\ref{T:Maintheorem} now follows from Corollary~\ref{C:spines} applied
to the surface $(X,\psi)$ constructed in Section~\ref{S:HypS} and 
logarithmic staircases over 0 and $\infty$. The points $z_j$ should be chosen
so that $\psi(z_j)$ belongs to a regularly covered interval of the real line
between a critical and an asymptotic values. 

\hfill\qed(Theorem~\ref{T:Maintheorem})

\end{document}